\documentclass[12pt]{article}
\usepackage{latexsym}
\usepackage{epsfig,amssymb,euscript}
\usepackage{amsmath}
\numberwithin{equation}{section}  

\newsavebox{\ns}
\newsavebox{\dbrane}
\newsavebox{\dbshort}

\def\be{\begin{equation}}
\def\ee{\end{equation}}
\def\bea{\begin{eqnarray}}
\def\eea{\end{eqnarray}}

\newcommand{\nn}{\nonumber}

\def\Dslash{\,\,{\raise.15ex\hbox{/}\mkern-12mu D}}
\def\Dbarslash{\,\,{\raise.15ex\hbox{/}\mkern-12mu {\bar D}}}
\def\delslash{\,\,{\raise.15ex\hbox{/}\mkern-9mu \partial}}
\def\delbarslash{\,\,{\raise.15ex\hbox{/}\mkern-9mu {\bar\partial}}}
\def\pslash{\,\,{\raise.15ex\hbox{/}\mkern-9mu p}}
\def\calDslash{\,\,{\raise.15ex\hbox{/}\mkern-12mu {\cal D}}}

\newcommand\R{\mathbb{R}}
\newcommand\Z{\mathbb{Z}}

\newcommand\C{\mathbb{C}}
\newcommand\T{\mathbb{T}}
\newcommand\diff{\mathrm{d}}

\newcommand{\vol}{\mathrm{vol}}

\newtheorem{theorem}{Theorem}[section] 
 
\newtheorem{proposition}[theorem]{Proposition} 
\newtheorem{corollary}[theorem]{Corollary} 
\newtheorem{conjecture}[theorem]{Conjecture}
\newtheorem{problem}[theorem]{Problem}

\newenvironment{definition}[1][Definition]{\begin{trivlist} 
\item[\hskip \labelsep {\bfseries #1}]}{\end{trivlist}} 
\newenvironment{example}[1][Example]{\begin{trivlist} 
\item[\hskip \labelsep {\bfseries #1}]}{\end{trivlist}} 
\newenvironment{remark}[1][Remark]{\begin{trivlist} 
\item[\hskip \labelsep {\bfseries #1}]}{\end{trivlist}} 
\newcommand{\qed}{\nobreak \ifvmode \relax \else \ifdim\lastskip<1.5em \hskip-\lastskip \hskip1.5em plus0em minus0.5em \fi \nobreak \vrule height0.75em width0.5em depth0.25em\fi}

\begin{document}
\begin{titlepage}
\begin{center}
\today

\vskip 2.5 cm 
{\Large \bf New Results in Sasaki-Einstein Geometry}\\

\vskip 1.5cm
{James Sparks}\\
\vskip 1 cm

Department of Mathematics, Harvard University \\
One Oxford Street, Cambridge, MA 02138, U.S.A.\\
\vskip 0.3cm
and
\vskip 0.3cm
Jefferson Physical Laboratory, Harvard University \\
Cambridge, MA 02138, U.S.A.\\

\vskip 1cm

\end{center}

\begin{abstract}
\noindent
This article is a summary of some of the author's work on Sasaki-Einstein geometry. A rather general conjecture in string theory known as the AdS/CFT 
correspondence relates Sasaki-Einstein geometry, in low dimensions, to 
superconformal field theory; properties of the latter 
are therefore reflected in the former, and vice versa. 
Despite this physical motivation, 
many recent results are of independent geometrical 
interest, and are described here in purely mathematical terms: explicit 
constructions of infinite families of both quasi-regular and irregular Sasaki-Einstein metrics; toric Sasakian geometry; 
an extremal problem that determines the Reeb vector field for, and hence also the volume of, a Sasaki-Einstein manifold; 
and finally, obstructions to the existence of Sasaki-Einstein metrics. Some of these results also provide new insights 
into K\"ahler geometry, and in particular new obstructions to the existence of K\"ahler-Einstein metrics on Fano orbifolds.

\vskip 1cm
\noindent {\it Invited contribution to the proceedings of the conference ``Riemannian Topology: Geometric Structures on Manifolds,'' 
Albuquerque, New Mexico, in celebration of Charles P. Boyer's 65$^{\mathrm{th}}$ birthday.}

\end{abstract}

\end{titlepage}
\pagestyle{plain}
\setcounter{page}{1}
\newcounter{bean}
\baselineskip18pt

\tableofcontents


\section{Introduction}

Sasaki-Einstein geometry 
is the odd-dimensional cousin of K\"ahler-Einstein 
geometry. In fact the latter, for positive Ricci curvature, is strictly
contained in the former; Sasaki-Einstein geometry is thus a 
generalisation of K\"ahler-Einstein geometry. 
The author's initial interest in this subject stemmed 
from a rather general conjecture in string theory known as 
the AdS/CFT correspondence \cite{Maldacena}. This is probably the 
most important conceptual development in theoretical physics 
in recent years. AdS/CFT conjecturally relates quantum 
gravity, in certain backgrounds, to ordinary quantum 
field theory without gravity. Moreover, the relation between the
two theories is \emph{holographic}: the quantum field theory 
resides on the \emph{boundary} of the region in which gravity propagates. 

In a particular setting the AdS/CFT correspondence
relates Sasaki-Einstein geometry, 
in dimensions five and seven, to superconformal field theory, in 
dimensions four and three, respectively. Superconformal field
theories are very special types of quantum field theories: 
they possess superconformal symmetry, and hence in particular 
conformal symmetry. The five-dimensional case of this correspondence
is currently understood best. One considers a ten-dimensional 
Riemannian product $(B\times L, g_B+g_L)$, where $(L,g_L)$ is a 
Sasaki-Einstein five-manifold and $(B,g_B)$ is five-dimensional 
\emph{hyperbolic space}. We may present $B$ as the open unit ball
$B = \{x\in \mathbb{R}^5\mid \|x\|<1\}$ with metric
\bea
g_B = \frac{4\sum_{i=1}^5 \diff x_i\otimes \diff x_i}{(1-\|x\|^2)^2}~.\eea
Here $\|x\|$ denotes the Euclidean norm of $x\in\mathbb{R}^5$. 
We may naturally compactify $(B,g_B)$, in the sense of Penrose, by adding a conformal boundary 
at infinity. One thus considers the closed unit ball 
$\bar{B}=\{x\in \mathbb{R}^5\mid \|x\|\leq 1\}$ equipped with the 
metric
\bea
g_{\bar{B}} = f^2 g_B\eea
where $f$ is a smooth function on $\bar{B}$ which is positive 
on $B$ and has a simple zero on $\partial \bar{B}=S^4$. For example, 
$f=1-\|x\|^2$ induces the standard metric on $S^4$, but there is no 
natural choice for $f$. Thus $S^4$ inherits only a 
\emph{conformal structure}. The isometric action of 
$SO(1,5)$ on $(B,g_B)$ extends to the action 
of the conformal group $SO(1,5)$ on the four-sphere. 
The AdS/CFT correspondence\footnote{The Lorentzian version of hyperbolic 
space is known as anti-de Sitter spacetime (AdS). The acronym CFT stands 
for conformal field theory.} 
conjectures that type IIB string theory,
which is supposed to be a theory of quantum gravity, propagating on 
$(B\times L,g_B+g_L)$ is equivalent to a four-dimensional 
superconformal field theory that resides on the boundary four-sphere. 
Indeed, the ten-dimensional manifold is a \emph{supersymmetric} 
solution to type IIB supergravity; in differentio-geometric terms, 
this means that there exists a solution to a 
certain Killing spinor equation. 
The AdS/CFT correspondence thus in particular implies a correspondence 
between Sasaki-Einstein manifolds in dimension five and superconformal 
field theories in four dimensions: for each Sasaki-Einstein 
five-manifold $(L,g_L)$ we obtain a different superconformal field 
theory. The AdS/CFT correspondence then naturally maps geometric 
properties of $(L,g_L)$ to properties of the dual superconformal field theory.

I should immediately emphasize,
however, that this article is aimed at geometers, rather than theoretical physicists. 
Unfortunately, explaining AdS/CFT to a mathematical audience is 
beyond the scope of the present paper. Instead I shall focus mainly
on the new geometrical results obtained by the author. 
The paper is based on a talk given at the conference 
``Riemannian Topology: Geometric Structures on Manifolds,'' in 
Albuquerque, New Mexico.

The outline of the paper is as follows. Section \ref{sasaki} contains
a brief review of Sasakian geometry, 
in the language of K\"ahler cones. Section \ref{explicit} 
summarises the properties of several infinite families of 
Sasaki-Einstein manifolds that were constructed in \cite{paper1, paper2, 
paper3, CLPP, CLPP2, MS2}, focusing on the (most physically interesting)
case of dimension five. Section \ref{toric} reviews 
toric Sasakian geometry, as developed in \cite{MSY}. 
Section \ref{Zmin} is a brief account of an extremal problem 
that determines the Reeb vector field for a Sasaki-Einstein metric 
\cite{MSY,MSY2}. This is understood completely for toric 
Sasaki-Einstein manifolds, whereas the general case currently contains some 
technical gaps. 
Finally, section \ref{obstructions} reviews some
obstructions to the existence of Sasaki-Einstein metrics \cite{GMSY}. 
A far more detailed account of Sasakian geometry, together with many 
other beautiful results in Sasaki-Einstein geometry not described here, 
may be found in the book \cite{BGbook}.


\section{Sasakian geometry}\label{sasaki}

While K\"ahler geometry \cite{kahler} has been studied intensively for more than seventy years, Sasakian geometry \cite{sasaki}
has, in contrast, received relatively little attention. Sasakian geometry 
was originally defined in terms of metric-contact geometry, but this does not really 
emphasize its relation to K\"ahler geometry. The following is a good
\begin{definition} A compact Riemannian manifold $(L,g_L)$ is \emph{Sasakian} if and only 
if its metric cone $(X_0=\mathbb{R}_+\times L, g=\diff r^2 + r^2g_L)$ is 
a \emph{K\"ahler cone}.\end{definition}
Here $r\in (0,\infty)$ may be regarded as a coordinate on the positive real line $\mathbb{R}_+$. 
The reason for the subscript on $X_0$ will become apparent later. Note that $(L,g_L)$ is 
isometrically embedded $\iota:L\rightarrow X_0$ into the cone with image $\{r=1\}$. By definition, $(X_0,g)$ is a non-compact  
K\"ahler manifold, with K\"ahler form 
\bea
\omega = \frac{1}{4}\diff \diff^c r^2~.\eea
Here $\diff^c = \mathcal{J}\circ\diff=i(\bar{\partial}-\partial)$, as usual, with $\mathcal{J}$ the complex structure tensor on $X_0$. 
The square of the radial function $r^2$ thus serves as a global K\"ahler potential on the 
cone. It is not difficult to verify that the homothetic vector field $r\partial/\partial r$ 
is holomorphic, and that 
\bea
\xi = \mathcal{J}\left(r\frac{\partial}{\partial r}\right)\eea
is holomorphic and also a \emph{Killing vector field}: $\mathcal{L}_{\xi}g=0$. 
$\xi$ is known as the \emph{Reeb vector field}. It is tangent to the surfaces 
of constant $r$, and thus defines a vector field on $L$, which, in a 
standard abuse of notation, we also denote by $\xi$. Another important object is the 
one-form
\bea
\eta = \diff^c \log r~.\eea
This is homogeneous degree zero\footnote{that is, $\mathcal{L}_{r\partial/\partial r}\eta=0$.} 
under $r\partial/\partial r$ and pulls back via $\iota^*$ to a one-form, 
which we also denote by $\eta$, on $L$. In fact $\eta$ is precisely a \emph{contact one-form} 
on $L$; that is, $\eta\wedge(\diff \eta)^{n-1}$ is a volume form on $L$. Here $n=\dim_{\mathbb{C}}X_0$,
or equivalently $2n-1=\dim_{\mathbb{R}} L$. The pair $\{\eta,\xi\}$ then satisfy, either on
the cone or on $L$, the relations 
\bea
\eta(\xi)=1,\qquad \diff \eta (\xi,\cdot) =0~.\eea
This is the usual 
definition of the Reeb vector field in contact geometry. Indeed, the open cone $(X_0,\omega)$ 
is the \emph{symplectization} of the contact manifold $(L,\eta)$, where one 
regards $\omega=\tfrac{1}{2}\diff (r^2\eta)$ as a symplectic form on $X_0$.

The square norm of $\xi$, in the cone metric $g$, is $\|\xi\|_g^2=r^2$, and thus in particular 
$\xi$ is nowhere zero on $X_0$. It follows that the orbits of $\xi$ define 
a foliation of $L$. It turns out that the metric transverse to these orbits $g_T$ is 
also a K\"ahler metric. Thus Sasakian structures are sandwiched between two K\"ahler 
structures: the K\"ahler cone of complex dimension $n$, and the transverse K\"ahler 
structure of complex dimension $n-1$. 

Consideration of the orbits of $\xi$ leads to a global
classification of Sasakian structures. Suppose that all the orbits of 
$\xi$ close. This means that $\xi$ generates an isometric $U(1)$ action on 
$(L,g_L)$. Such Sasakian structures are called \emph{quasi-regular}. Since $\xi$ is nowhere zero, this action must be locally free: 
the isotropy subgroup at any point must be finite, and therefore isomorphic 
to a cyclic group $\mathbb{Z}_m\subset U(1)$. If the isotropy subgroups for 
all points are trivial then the $U(1)$ action is free, and the Sasakian 
structure is called \emph{regular}. We use the term \emph{strictly quasi-regular} 
for a quasi-regular Sasakian structure that is not regular.
In either case there is a quotient $V=L/U(1)$, which is generally an orbifold. The isotropy subgroups descend to 
the local orbifold structure groups in the quotient space $V=L/U(1)$; thus $V$ is a manifold
when the Sasakian structure is regular. The transverse K\"ahler metric descends to a K\"ahler metric 
on the quotient, so that $(V,g_V)$ is a K\"ahler manifold or orbifold. 
Indeed, $(V,g_V)$ may be regarded as the \emph{K\"ahler reduction} of the 
K\"ahler cone $(X_0,\omega)$ with respect to the $U(1)$ action, which is Hamiltonian 
with Hamiltonian function $\tfrac{1}{2}r^2$. 

If the orbits of $\xi$ do not all close, the Sasakian structure is said to be 
\emph{irregular}. The generic orbit is $\mathbb{R}$, and in this case 
one cannot take a meaningful quotient. The closure of the orbits of $\xi$ 
defines an abelian subgroup of the isometry group of $(L,g_L)$. Since $L$ is compact, 
the isometry group of $(L,g_L)$ is compact, and the closure of the orbits of 
$\xi$ therefore defines a torus $\mathbb{T}^s$, $s>1$, which acts 
isometrically on $(X_0,g)$ or $(L,g_L)$. Thus irregular Sasakian manifolds 
have at least a $\mathbb{T}^2$ isometry group.

The main focus of this article will be \emph{Sasaki-Einstein} manifolds. A simple 
calculation shows that\footnote{Notice the slight abuse of notation here: we 
are regarding all tensors in this equation as tensors on $X_0$.}
\bea
\mathrm{Ric}(g) = \mathrm{Ric}(g_L)-2(n-1)g_L = \mathrm{Ric}(g_T)-2ng_T\eea
where $\mathrm{Ric}(\cdot)$ denotes the Ricci tensor of a given metric. Thus the K\"ahler cone 
$(X_0,g)$ is Ricci-flat if and only if $(L,g_L)$ is Einstein with positive scalar
curvature $2(n-1)(2n-1)$, if and only if the transverse metric is Einstein with positive scalar curvature $4n(n-1)$. 
$(L,g_L)$ is then said to be a Sasaki-Einstein manifold. Notice that 
$(X_0,g)$ is then a Calabi-Yau cone, and in the quasi-regular case the 
circle quotient of $(L,g_L)$ is a K\"ahler-Einstein manifold or orbifold of positive Ricci curvature. 
The converse is also true: give a positively curved K\"ahler-Einstein manifold or orbifold $(V,g_V)$, there 
exists a Sasaki-Einstein metric\footnote{In the orbifold case this lifting 
may or may not be an orbifold. If $\{\Gamma_{\alpha}\}$ denote the 
local orbifold structure groups of $V$, then the data that defines 
an orbibundle over $V$ with structure group $G$ includes elements of $\mathrm{Hom}(\Gamma_{\alpha},G)$ 
for each $\alpha$, subject to certain gluing conditions. A moment's thought shows 
that the total space of a $G$ principal orbibundle over $V$ is smooth 
if and only if all these maps are \emph{injective}.} 
on the total space of a $U(1)$ principal (orbi-)bundle over $V$. 
This was proven in general by Boyer and Galicki in \cite{boyer}. 


\section{Explicit constructions of Sasaki-Einstein manifolds}\label{explicit}

Explicit examples of Sasaki-Einstein manifolds were, until recently, quite 
rare. In dimension five, the only simply-connected\footnote{Positively 
curved Einstein manifolds have finite fundamental group \cite{myers}.} examples 
that were known in explicit form were the round sphere and a certain 
homogeneous metric on $S^2\times S^3$ \cite{tanno}, known as $T^{1,1}$ in the physics 
literature. These are both regular, being circle bundles over 
$\mathbb{CP}^2$ and $\mathbb{CP}^2\times\mathbb{CP}^1$ with their 
standard K\"ahler-Einstein metrics, respectively. In fact regular 
Sasaki-Einstein manifolds are classified~\cite{fried}. This follows since 
smooth K\"ahler-Einstein surfaces with positive Ricci
curvature have been classified by Tian and Yau
\cite{tian,tianyau}. The result is that the base may be
taken to be a del Pezzo
surface obtained by blowing up $\mathbb{CP}^2$ at $k$ generic
points with $3\le k\le8$; although proven to exist, the K\"ahler-Einstein 
metrics on these del Pezzo surfaces are not known explicitly. More 
recently, Boyer, Galicki and collaborators have produced vast numbers of 
quasi-regular Sasaki-Einstein metrics using existence results of Koll\'ar for K\"ahler-Einstein metrics on 
Fano orbifolds, together with the $U(1)$ lifting 
mentioned at the end of the previous section. For a review of their work, see 
\cite{BGreview,BGbook}.

Until 2004 no explicit examples of non-trivial strictly quasi-regular Sasaki-Einstein manifolds 
were known, and it was not known whether or not \emph{irregular} Sasaki-Einstein 
manifolds even existed. In fact Cheeger and Tian conjectured in \cite{cheeger} 
that they did not exist. The following theorem disproves this conjecture:
\begin{theorem}\label{Ypq}(\cite{paper2}) There exist a countably infinite number of 
Sasaki-Einstein metrics $Y^{p,q}$ on $S^2\times S^3$, labelled naturally by $p,q\in \mathbb{N}$ where 
$\mathrm{gcd}(p,q)=1$, $q<p$. $Y^{p,q}$ is quasi-regular if and only if 
$4p^2-3q^2$ is the square of a natural number, otherwise it is irregular. 
In particular, there are infinitely many strictly quasi-regular and irregular 
Sasaki-Einstein metrics on $S^2\times S^3$.
\end{theorem}
These metrics were constructed explicitly in \cite{paper2}, based on 
supergravity constructions by the same authors in \cite{paper1}. 
The metrics are cohomogeneity one, meaning that the generic orbit under 
the action of the isometry group has real codimension one. 
The Lie algebra of this group is $\mathtt{su}(2)\times\mathtt{u}(1)\times\mathtt{u}(1)$.
The volumes of the metrics are given by 
the formula
\bea
\frac{\vol[Y^{p,q}]}{\pi^3} = \frac{q^2(2p+\sqrt{4p^2-3q^2})}{3p^2(3q^2-2p^2+\sqrt{4p^2-3q^2})}~.\eea
The result that $Y^{p,q}$ is diffeomorphic to $S^2\times S^3$ follows 
from Smale's classification of 5-manifolds \cite{smale}. Interestingly, the cone $X_0$ corresponding to $Y^{2,1}$ is 
the open complex cone over the first del Pezzo surface \cite{MS}. Note that the first del Pezzo surface was missing from the list of Tian and Yau; 
it cannot admit a K\"ahler-Einstein metric since its Futaki invariant is non-zero. The Ricci-flat K\"ahler cone 
metric is in fact \emph{irregular} for $Y^{2,1}$. We shall return to this in section \ref{obstructions}.
Recently, Conti \cite{conti} has classified cohomogeneity one Sasaki-Einstein 
five-manifolds: they are precisely the set $\{Y^{p,q}\}$. The construction of the above
metrics also easily extends to higher dimensions \cite{paper3,general1,general2}. This leads to
the following 
\begin{corollary} There exist countably infinitely many strictly quasi-regular and irregular Sasaki-Einstein structures 
in every (odd) dimension greater than 3.\end{corollary}
This corollary should be contrasted with several other results. It is known that for fixed dimension there 
are finitely many (deformation classes of) Fano manifolds \cite{kollar}; thus there are only finitely many 
positively curved K\"ahler-Einstein structures in each dimension, and hence finitely many
\emph{regular} Sasaki-Einstein structures in each odd dimension. On the other hand, 
these may occur in continuous families. This is already true for del Pezzo surfaces 
with $k\geq 5$ blow-ups, which have a complex structure moduli space of complex dimension 
$2(k-4)$. The quasi-regular existence results of Boyer and Galicki also produce examples 
with, sometimes quite large, moduli spaces \cite{BGreview}. It is currently unknown whether 
or not there exist continuous families of \emph{irregular} Sasaki-Einstein structures.

Perhaps surprisingly, there also exist explicit cohomogeneity two Sasaki-Einstein five-manifolds. 
The following subsumes Theorem \ref{Ypq}:
\begin{theorem}\label{Labc}(\cite{CLPP},\cite{CLPP2},\cite{MS2})  \ \ There exist a countably infinite number of 
Sasaki-Einstein metrics $L^{a,b,c}$ on $S^2\times S^3$, labelled naturally by $a,b,c \in \mathbb{N}$ where
$a\leq b$, $c\leq b$, $d=a+b-c$, $\mathrm{gcd}(a,b,c,d)=1$, $\mathrm{gcd}(\{a,b\},\{c,d\})=1$. 
Here the latter means that each of the pair $\{a,b\}$ must be coprime to 
each of $\{c,d\}$. Moreover, $L^{p-q,p+q,p}=Y^{p,q}$. The metrics are generically 
cohomogeneity two, generically irregular, and generically have volumes that 
are the product of quartic irrational numbers with $\pi^3$.
\end{theorem}
The condition under which the metrics are (strictly) quasi-regular is not simple to determine in general. 
The quartic equation with integer coefficients that is 
satisfied by $\vol[L^{a,b,c}]/\pi^3$ is written down explicitly in \cite{CLPP}. 
For integers $(a,b,c)$ not satisfying some of the coprime conditions 
one obtains Sasaki-Einstein orbifolds. 

We conclude this section with a comment on how the metrics in Theorem \ref{Labc} 
were constructed. In \cite{CLPP,CLPP2} the local form of the metrics was found 
by writing down the Riemannian forms of known black hole metrics, and then taking 
a certain ``BPS'' limit. The initial family of metrics are local Einstein metrics, and in the limit 
one obtains a local family of Sasaki-Einstein metrics. One then determines when these local metrics
extend to complete metrics on compact manifolds, and this is where the integers $(a,b,c)$ enter.
The same metrics were independently discovered in a slightly different manner in \cite{MS2}. 
In the latter reference the local K\"ahler-Einstein metrics in dimension four are constructed first. 
It turns out that these are precisely the \emph{orthotoric} K\"ahler-Einstein metrics 
in \cite{apostolov}. The construction again easily extends to 
higher dimensions \cite{CLPP,CLPP2}; the metrics are generically 
cohomogeneity $n-1$.


\section{Toric Sasakian geometry}\label{toric}

In this section we summarise some of the results in \cite{MSY} on toric Sasakian 
geometry. This probably warrants a 
\begin{definition} A Sasakian manifold $(L,g_L)$ is said to be \emph{toric} 
if there exists an effective, holomorphic and Hamiltonian action 
of the torus $\T^n$ on the corresponding K\"ahler cone $(X_0,g)$. 
The Reeb vector field $\xi$ is assumed to lie in the Lie algebra 
of the torus $\xi\in \mathtt{t}_n$. \end{definition}
The Hamiltonian condition means that there exists a $\T^n$-invariant moment map 
\bea
\mu:X_0\rightarrow \mathtt{t}^*_n~.\eea
The condition on the Reeb vector field implies
that the image is a strictly convex rational polyhedral cone \cite{FT, lerman}. 
Symplectic toric cones with Reeb vector fields not satisfying this condition form a short list and
have been classified in \cite{lerman}. The main result of this section is 
Proposition \ref{class} which describes, in a rather explicit form, the 
space of toric Sasakian metrics on the link of a fixed affine toric  
variety $X$. Any toric Sasakian manifold is of this form, with 
the open K\"ahler cone $X_0=X\setminus \{p\}$ being the smooth part of 
$X$, with $p$ an isolated singular point.

We begin by fixing a \emph{strictly convex rational polyhedral cone} $\mathcal{C}^*$ 
in $\mathbb{R}^n$, where the latter is regarded as the dual Lie algebra of a torus $\mathtt{t}_n^*\cong\mathbb{R}^n$ 
with a particular choice of basis:
\bea
\mathcal{C}^* = \{y\in \mathtt{t}_n^*\mid \langle y,v_a\rangle \geq 0, \forall a=1,\ldots,d\}~.\eea
The strictly convex condition means that $\mathcal{C}^*$ is a cone over a convex polytope of 
dimension $n-1$. It follows that necessarily $n\leq d\in \mathbb{N}$. The rational 
condition on $\mathcal{C}^*$ means that the vectors $v_a\in\mathtt{t}_n\cong\mathbb{R}^n$ are rational. 
In particular, one can normalise the $v_a$ so that they are primitive vectors in $\mathbb{Z}^n \cong \ker \{\exp: 
\mathtt{t}_n\rightarrow \mathbb{T}^n\}$. The $v_a$ are thus the inward-pointing primitive normal 
vectors to the bounding hyperplanes of the polyhedral cone $\mathcal{C}^*$. We may alternatively 
define $\mathcal{C}^*$ in terms of its generating vectors $\{u_{\alpha}\in\mathbb{Z}^n\}$:
\bea
\mathcal{C}^* = \left\{\sum_{\alpha} \lambda_{\alpha} u_{\alpha}\mid \lambda_{\alpha}\geq 0\right\}~.\eea
The primitive vectors $u_{\alpha}$ generate the one-dimensional faces, or \emph{rays}, of the polyhedral cone 
$\mathcal{C}^*$.

Define the linear map
\bea A:&& \mathbb{R}^d\rightarrow \mathbb{R}^n\nonumber\\
&& e_a\mapsto v_a\eea
where $\{e_a\}$ denotes the standard orthonormal basis of $\mathbb{R}^d$. Let $\Lambda\subset\mathbb{Z}^n$ 
denote the lattice spanned by $\{v_a\}$ over $\mathbb{Z}$. This is of maximal rank, since $\mathcal{C}^*$ 
is strictly convex. There is an induced map of tori
\bea
\mathbb{T}^d\cong \mathbb{R}^d/2\pi\mathbb{Z}^d\rightarrow \mathbb{R}^n/2\pi\mathbb{Z}^n\cong\mathbb{T}^n\eea
where the kernel is a compact abelian group 
$\mathcal{A}$, with $\pi_0(\mathcal{A})\cong \Gamma\cong \mathbb{Z}^n/\Lambda$. 

Using this data we may construct the following K\"ahler quotient:
\bea\label{quotient}
X = \mathbb{C}^d//\mathcal{A}~.\eea
Here we equip $\mathbb{C}^d$ with its standard flat K\"ahler structure $\omega_{\mathrm{flat}}$. $\mathcal{A}\subset 
\mathbb{T}^d$ acts holomorphically and Hamiltonianly on $(\mathbb{C}^d,\omega_{\mathrm{flat}})$. 
We then take the K\"ahler quotient (\ref{quotient}) at level zero. The origin of $\C^d$ 
projects to a singular point in $X$, and the induced K\"ahler metric 
$g_{\mathrm{can}}$ on its complement $X_0$ is a cone. Moreover, the 
quotient torus $\mathbb{T}^d/\mathcal{A}\cong\mathbb{T}^n$ acts holomorphically 
and Hamiltonianly on $(X_0,\omega_{\mathrm{can}})$, with moment map 
\bea
\mu:X_0\rightarrow\mathtt{t}_n^*; \qquad \mu(X)=\mathcal{C}^*~.
\eea
The quotient (\ref{quotient}) may be written explicitly as follows. One computes 
a primitive basis for the kernel of $A$ over $\Z$ by finding all 
solutions to 
\bea
\sum_a Q_{I}^a  v_a = 0\eea
with $Q_{I}^a\in \Z$, and such that for each $I$ the $\{Q_I^a\mid a=1,\ldots,d\}$ have no 
common factor. The number of solutions, which are indexed by $I$, is 
$d-n$ since $A$ is surjective; this latter fact again
follows since $\mathcal{C}^*$ is strictly convex. One then has 
\bea\label{reduction}
X = \mathcal{K}/\mathcal{A}
\equiv \C^d//\mathcal{A}\eea
with
\bea
\mathcal{K} \equiv \left\{(Z_1,\ldots,Z_d)\in\C^d\mid\sum_{a}Q_{I}^a |Z_a|^2 = 0\right\}\subset \C^d\eea
where $Z_a$ denote standard complex coordinates on $\C^d$ and the charge matrix 
$Q_{I}^a$ specifies the torus embedding $\T^{d-n}\subset \T^d$. 

It is a standard fact that the space $X$ is an affine toric variety; that is, $X$ is an affine variety 
equipped with an effective holomorphic action of the \emph{complex} torus $\mathbb{T}^n_{\mathbb{C}}\cong (\mathbb{C}^*)^n$ 
which has a dense open orbit. Let 
\bea
\mathcal{C}=\{\xi\in\mathtt{t}_n\mid \langle \xi,y\rangle\geq 0, \forall y\in\mathcal{C}^*\}~.\eea
This is the \emph{dual cone} to $\mathcal{C}^*$, which is also 
a convex rational polyhedral cone by Farkas' Theorem. In the algebro-geometric language, 
the cone $\mathcal{C}$ is precisely the \emph{fan} for the affine toric variety $X$. 
We have $X_0=\mathbb{R}_+\times L$ with $L$ compact. If one begins with a general 
strictly convex rational polyhedral cone $\mathcal{C}^*$, the link $L$ will be an orbifold; in order that $L$ be a smooth manifold one 
requires the moment polyhedral cone $\mathcal{C}^*$ to be \emph{good} \cite{lerman}. This puts certain 
additional constraints on the vectors $v_a$; the reader is
referred to \cite{MSY} for the details. Note that 
$L$ inherits a canonical Sasakian metric from the K\"ahler quotient metric $g_{\mathrm{can}}$ on $X_0$. 

Let $\partial/\partial\phi_i$, $i=1,\ldots,n$, be a 
basis\footnote{By a standard abuse of notation we identify vector fields on 
$X_0$ with corresponding elements of the Lie algebra.} for $\mathtt{t}_n$, where 
$\phi_i\in [0,2\pi)$ are coordinates on the real torus $\mathbb{T}^n$. Then we have the following
\begin{proposition}\label{class}(\cite{MSY}) The space of 
toric K\"ahler cone metrics on the smooth part of an affine toric variety $X_0$ is a product
\bea
\mathcal{C}_{\mathrm{int}}\times \mathcal{H}^1(\mathcal{C}^*)\nonumber\eea
where $\xi\in\mathcal{C}_{\mathrm{int}}\subset\mathtt{t}_n$ labels the Reeb vector field, with 
$\mathcal{C}_{\mathrm{int}}$ the open interior of $\mathcal{C}$, and 
$\mathcal{H}^1(\mathcal{C}^*)$ denotes the space of homogeneous degree one functions on $\mathcal{C}^*$ that are 
smooth up to the boundary (together with the convexity condition below). 

Explicitly, on the 
dense open image of $\mathbb{T}^n_{\mathbb{C}}$ we have
\bea
g = G_{ij}\diff y^i\diff y^j + G^{ij}\diff\phi_i\diff\phi_j\eea
where 
\bea
G_{ij}=\frac{\partial^2G}{\partial y^i\partial y^j}\eea
with matrix inverse $G^{ij}$, and the function
\bea
G(y) = G_{\mathrm{can}}(y) + G_{\xi}(y) + h(y)\eea
is required to be strictly convex with $h(y)\in \mathcal{H}^1(\mathcal{C}^*)$ and
\begin{eqnarray}
G_{\mathrm{can}}(y) & =& \frac{1}{2}\sum_{a=1}^d \langle y,v_a\rangle\log \langle y,v_a\rangle\nonumber\\
G_{\xi}(y) &= &\frac{1}{2}\langle\xi,y\rangle\log \langle\xi,y\rangle - \frac{1}{2}\left(\sum_{a=1}^d 
\langle v_a,y\rangle\right)\log \left(\sum_{a=1}^d 
\langle v_a,y\rangle\right)~.\nonumber\end{eqnarray}
\end{proposition}
In particular, the canonical metric on $X_0$, induced from K\"ahler reduction of 
the flat metric on $\mathbb{C}^d$, is given by setting $G(y)=G_{\mathrm{can}}(y)$. 
This function has a certain singular behaviour at the boundary $\partial\mathcal{C}^*$ 
of the polyhedral cone; this is required precisely so that the metric compactifies 
to a smooth metric on $X_0$. 

The space of Reeb vector fields is the interior of 
the cone $\mathcal{C}$. One can show \cite{MSY2} that for $\xi\in\partial\mathcal{C}$ 
the vector field $\xi$ must vanish somewhere on $X_0$. Specifically, the bounding facets 
of $\mathcal{C}$ correspond to the generating rays of $\mathcal{C}^*$ under the duality map between 
cones; $\xi$ being in a bounding facet of $\mathcal{C}$ implies that the corresponding vector field
then vanishes on the inverse image, under the moment map, of the dual generating ray of $\mathcal{C}^*$. 
However, since the Reeb vector field is nowhere vanishing, we see that the boundary of 
$\mathcal{C}$ is a singular limit of Sasakian metrics on $X_0$.

Fixing a particular choice of K\"ahler cone metric on $X_0$, the image of $L=\{r=1\}$
under the moment map is
\bea
\mu(L) = \left\{y\in\mathcal{C}^*\mid \langle y,\xi \rangle =\tfrac{1}{2}\right\}\equiv H({\xi})~.\eea
The hyperplane $\langle y,\xi\rangle = \tfrac{1}{2}$ is called the \emph{characteristic 
hyperplane} \cite{boyertoric}. This intersects the moment cone $\mathcal{C}^*$ to form a 
compact $n$-dimensional polytope $\Delta({\xi})=\mu(\{r\leq 1\})$, bounded by $\partial\mathcal{C}^*$ 
and the compact $(n-1)$-dimensional polytope $H({\xi})$. In particular, the image $H(\xi)$ of $L$ under
the moment map depends only on the Reeb vector field $\xi$, and not on the choice of 
homogeneous degree one function $h$ in Proposition \ref{class}. Moreover, the 
\emph{volume} of a toric Sasakian manifold $(L,g_L)$ is \cite{MSY}
\bea\label{euclid}
\vol[g_L] = 2n (2\pi)^n\vol[\Delta({\xi})]\eea
where $\vol[\Delta({\xi})]$ is the \emph{Euclidean} volume of $\Delta({\xi})$.

Finally in this section we introduce the notion of a \emph{Gorenstein singularity}:
\begin{definition} An analytic space $X$ with isolated singular point $p$ and smooth part $X\setminus\{p\}=X_0$
is said to be \emph{Gorenstein} if there exists a smooth nowhere zero holomorphic $(n,0)$-form $\Omega$ 
on $X_0$.\end{definition}
We shall refer to $\Omega$ as a \emph{holomorphic volume form}. $X$ being Gorenstein 
is a necessary condition for $X_0$ to admit a Ricci-flat K\"ahler metric, and hence for
the link $L$ to admit a Sasaki-Einstein metric. Indeed, the Ricci-form $\rho=\mathrm{Ric}(\mathcal{J}\cdot,\cdot)$ is a curvature two-form for the holomorphic line bundle $\Lambda^{n,0}$. The Ricci-flat 
K\"ahler condition implies
\bea\label{MA}
\frac{i^n}{2^n}(-1)^{n(n-1)/2}\Omega\wedge\bar{\Omega} = \frac{1}{n!}\omega^n~.\eea
For affine toric varieties, it is again well-known that $X$ being
Gorenstein is equivalent to the existence of a basis for the torus $\mathbb{T}^n$ 
for which $v_a=(1,w_a)$ for each $a=1,\ldots,d$, and $w_a\in\mathbb{Z}^{n-1}$. 

\begin{example}(\cite{MS,MS2}) From Theorems \ref{Ypq} and \ref{Labc} one sees
that $L^{a,b,c}$, which contain $Y^{p,q}$ as a subset, have a holomorphic 
Hamiltonian action of $\mathbb{T}^3$ on the corresponding K\"ahler cones and are thus toric Sasaki-Einstein manifolds. 
One finds 
that the image of the cone under the moment map is always a four--sided polyhedral cone ($d=4$) in $\mathbb{R}^3$. The charge 
matrix $Q$ is 
\bea
Q=(a,b,-c,-a-b+c)~.\eea
The Gorenstein condition is reflected by the fact that the sum of the components of $Q$ is zero. 
In particular, for $Y^{p,q}$ (which is $a=p-q$, $b=p+q$, $c=p$) we have
\bea
v_1=[1,0,0],\quad v_2=[1,1,0],\quad v_3=[1,p,p], \quad v_4=[1,p-q-1,p-q]~.\eea
It is relatively straightforward to see that the affine toric Gorenstein singularities for
$L^{a,b,c}$ are the most general such that are generated by four rays.
\end{example}


\section{A variational problem for the Reeb vector field}\label{Zmin}

In this section we consider the following problem: given a
Gorenstein singularity $(X,\Omega)$, with isolated singular 
point and smooth set $X_0=\mathbb{R}_+\times L$, what is the Reeb vector 
field for a Ricci-flat K\"ahler cone metric on $X_0$, assuming it exists? We shall go 
quite a long way in answering this question, and give a complete solution 
for affine toric varieties; in general more work still remains to be done.

The strategy is to set up a variational problem on a space of Sasakian metrics on $L$,
 or equivalently a space of K\"ahler cone metrics on $X_0$. To this end, we 
suppose that $X_0$ is equipped with an effective holomorphic action of the torus $\mathbb{T}^s$ 
for some $s$; this is clearly necessary. Indeed, for irregular Sasakian metrics 
one requires $s>1$, as commented in section \ref{sasaki}. We then assume we are given a space of K\"ahler cone metrics $\mathcal{S}(X_0)$
on $X_0$ such that:
\begin{itemize}
\item The torus $\mathbb{T}^s$ acts Hamiltonianly on each metric $g\in \mathcal{S}(X_0)$.
\item The Reeb vector field for each metric lies in the Lie algebra $\mathtt{t}_s$ of $\mathbb{T}^s$.
\end{itemize}
We shall  
continue to denote K\"ahler cone metrics by $g$, the corresponding 
Sasakian metric by $g_L$, and regard either as elements of $\mathcal{S}(X_0)$.
The second condition above ensures that the torus action is of \emph{Reeb type} \cite{FT,boyertoric}. 
We then have the following
\begin{proposition}\label{reeby}(\cite{MSY2}) The \emph{volume} of the link $(L,g_L)$, 
as a functional on the space $\mathcal{S}(X_0)$, depends only on the Reeb vector field $\xi$ 
for the Sasakian metric $g_L\in \mathcal{S}(X_0)$.\end{proposition}
It follows that $\vol$ may be regarded as a function on the space of Reeb vector fields: 
\bea
\vol:\mathcal{R}(X_0)\rightarrow \mathbb{R}_+\eea
where
\bea
\mathcal{R}(X_0)=\left\{\xi\in\mathtt{t}_s\mid \xi = \mathrm{Reeb} \ \mathrm{vector} \ \mathrm{field} \ \mathrm{for} \ \mathrm{some} \ 
g_L\in\mathcal{S}(X_0)\right\}~.\eea
The first and second derivatives are given by
\begin{proposition}\label{der}
\bea\label{dev1}
\diff\vol(Y)&=& -n\int_L \eta(Y)\diff\mu\\\label{dev2}
\diff^2\vol(Y,Z)& =&n(n+1)\int_L \eta(Y)\eta(Z)\diff\mu~.\eea
Here $Y,Z$ are holomorphic Killing vector fields in $\mathtt{t}_s$, 
$\eta$ is the contact one-form for the Sasakian metric, and $\diff\mu$ is the 
Riemannian measure on $(L,g_L)$.\end{proposition}
Note that for toric Sasakian metrics (for which the torus $\mathbb{T}^s$ has maximal dimension: $s=n$)
we already noted Proposition \ref{reeby} in the previous section -- 
see equation (\ref{euclid}). Note also that 
(\ref{dev2}) shows that $\vol$ is a strictly convex function of $\xi$. 

Proposition \ref{reeby} is proven roughly as follows. Suppose one has two 
K\"ahler cone metrics on $X_0$ with the same homothetic vector field $r\partial/\partial r$. 
It is then straightforward to show that the K\"ahler potentials differ by 
a multiplicative factor $\exp\varphi$, where $\varphi$ is a basic homogeneous degree zero 
function: 
$\mathcal{L}_{\xi}\varphi=0=\mathcal{L}_{r\partial/\partial r} \varphi$, where 
recall that $\xi=\mathcal{J}(r\partial/\partial r)$. One then shows
that the volume is independent of $\varphi$.

This is precisely 
analogous to the situation in K\"ahler geometry where one fixes a \emph{K\"ahler class}. 
Suppose that $(M,\omega)$ is a \emph{compact} K\"ahler manifold with K\"ahler class $[\omega]\in H^{1,1}(M)$. 
Then any other K\"ahler metric on $M$ in the same K\"ahler class is given by 
$\omega + i\partial\bar{\partial}\varphi$ for some smooth real function $\varphi$. 
The volume of $(M,\omega)$ clearly depends only on $[\omega]$. 

Indeed, one can push the analogy further. A choice of Reeb vector field 
on $X_0$ should be regarded as a choice of \emph{polarisation}\footnote{This 
terminology was introduced in \cite{BGS}.}. The space of quasi-regular 
Reeb vector fields is dense in the space of all Reeb vector fields: 
quasi-regular Reeb vector fields correspond to rational vectors in the Lie algebra 
$\mathtt{t}_s\cong\mathbb{R}^s$, and these are dense since the rationals are dense in the reals.  For $\xi$ quasi-regular, the 
$U(1)$ quotient is a K\"ahler orbifold $(V,\omega_V)$. Changing the polarisation 
$\xi$ thus changes the quotient $V$, in contrast to the K\"ahler setting in the last paragraph
where $M$ is fixed and the K\"ahler class changes. It is also straightforward
to show that the space of Reeb vector fields $\xi$ forms a \emph{cone}: 
if $\xi$ is a Reeb vector field, then $c\xi$ is also a Reeb vector field 
for another K\"ahler cone metric on $X_0$, for any constant $c>0$. Thus 
the space $\mathcal{R}(X_0)$ of Reeb polarisations forms a cone, analogous to the 
K\"ahler cone in K\"ahler geometry: we saw this explicitly 
in the previous section on toric Sasakian manifolds, where the space of 
Reeb vector fields is the interior $\mathcal{C}_{\mathrm{int}}$ of the 
polyhedral cone $\mathcal{C}$. 

We now suppose that $X$ is also Gorenstein. This is necessary for the 
existence of a Ricci-flat K\"ahler metric on $X_0$. We then introduce
the subspace $\mathcal{S}(X_0,\Omega)$ as the space of metrics 
in $\mathcal{S}(X_0)$ for which the Reeb vector field $\xi$ 
satisfies
\bea\label{charge}
\mathcal{L}_{\xi}\Omega = in\Omega~.\eea
Equivalently, $\Omega$ should be homogeneous degree $n$ under $r\partial/\partial r$. 
This is again clearly a necessary condition for a Ricci-flat K\"ahler cone metric, {\it cf} (\ref{MA}).

Recall now that Einstein metrics on $L$ are critical points of the \emph{Einstein-Hilbert
action}:
\bea
\mathcal{I}  : \mathrm{Metrics}(L) & \rightarrow &\mathbb{R}\nn\\
 g_L &\mapsto &\int_L \left[s(g_L)+2(n-1)(3-2n)\right]\diff\mu\eea
where $s(g_L)$ is the scalar curvature of $g_L$. We then have the following proposition:
\begin{proposition}(\cite{MSY2})
The Einstein-Hilbert action, as a functional on $\mathcal{S}(X_0)$,
depends only on the Reeb vector field $\xi$. It may thus be regarded as a function 
of $\xi$. Moreover, for Sasakian metrics $g_L\in\mathcal{S}(X_0,\Omega)$ we have
\bea
\mathcal{I}(g_L) = 4(n-1)\vol[g_L]~.\eea
\end{proposition}
Thus the Einstein-Hilbert action restricted to the space $\mathcal{S}(X_0,\Omega)$ 
is simply the volume functional, and depends only on the Reeb vector field $\xi$ of the metric. 
This suggests we introduce
\bea
\mathcal{R}(X_0,\Omega) = \left\{\xi\in\mathcal{R}(X_0)\mid \mathcal{L}_{\xi}\Omega= in\Omega\right\}~.\eea
Since Sasaki-Einstein metrics are critical points of $\mathcal{I}$, we see that 
the Reeb vector field for a Sasaki-Einstein metric is determined by a \emph{finite-dimensional 
extremal problem}, namely $\diff \mathcal{I}=0$, where $\mathcal{I}$ is interpreted as a 
function on $\mathcal{R}(X_0,\Omega)$. For toric varieties this is particularly simple:
\begin{theorem}\label{torictheorem} Let $X$ be an affine toric Gorenstein variety with 
fan (or Reeb polytope) $\mathcal{C}\subset\mathtt{t}_n$, and 
generating vectors of the form $v_a=(1,w_a)$. Then 
the Reeb vector field $\xi$ for a Ricci-flat K\"ahler cone metric 
on $X_0$ is uniquely determined as the critical point of the 
Euclidean volume of the polytope $\Delta({\xi})$
\bea
\vol[\Delta]:N_{\mathrm{int}}\rightarrow\mathbb{R}_+\eea
where $N$ is the $(n-1)$-dimensional polytope $N=\{\xi\in\mathcal{C}\mid 
\langle (1,0,\ldots,0), \xi \rangle =n\}$. 
\end{theorem}
Here we may take $\mathcal{R}(X_0)=\mathcal{C}_{\mathrm{int}}$ and $\mathcal{R}(X_0,\Omega)=N_{\mathrm{int}}$, following 
the classification of toric Sasakian metrics in section \ref{toric}. 
The first and second derivatives 
(\ref{dev1}), (\ref{dev2}) in Proposition \ref{der} may be written 
\cite{MSY}
\bea\label{toricdev1}
\frac{\partial\vol[\Delta]}{\partial\xi_i} & = &\frac{1}{2 \xi_k\xi_k}\int_{H({\xi})} 
y^i \ \diff \sigma\\ \label{toricdev2}
\frac{\partial^2\vol[\Delta]}{\partial\xi_i\partial \xi_j} &= &\frac{2(n+1)}{\xi_k\xi_k}\int_{H({\xi})} 
y^i y^j \ \diff\sigma~.\eea
Here $\diff\sigma$ is the standard measure induced on the $(n-1)$-polytope $H({\xi})\subset\mathcal{C}^*$. 
Uniqueness and existence of the critical point follows from a standard 
convexity argument: $\vol[\Delta]$ is a strictly convex (by (\ref{toricdev2})) positive function on 
the interior of a compact convex polytope 
$N$. Moreover, $\vol[\Delta]$ diverges to $+\infty$ at $\partial N$. It follows that $\vol[\Delta]$ 
must have precisely one critical point in the interior of $N$.

Given these results, it is natural to consider the general case, with $s\leq n$. 
We begin by recalling a classical result. Let $-\nabla^2_L$ be the scalar Laplacian 
on $(L,g_L)$, with spectrum $\{E_{\nu}\}_{\nu=0}^{\infty}$. Then we may define 
the heat kernel trace
\bea\label{heat}
\Theta(t) = \sum_{\nu=0}^{\infty} \exp(-t E_{\nu})\eea
where $t\in (0,\infty)$. There is a holomorphic analogue of this. Let $f$ be a holomorphic function on $X_0$ with 
\bea
\mathcal{L}_{\xi} f = \lambda if\eea
where $\R\ni \lambda>0$, 
and we refer to $\lambda$ as the charge of $f$ under 
$\xi$. Since $f$ is holomorphic, this immediately implies that
\bea
f = r^{\lambda}\tilde{f}\eea
where $\tilde{f}$ is homogeneous degree zero under $r\partial/\partial r$; that is, $\tilde{f}$ is the pull--back to $X$ 
of a function on the link $L$. 
Moreover, since $(X_0,g)$ is K\"ahler one can show that $f$ is \emph{harmonic}, and that 
\bea
-\nabla^2_L \tilde{f} = E\tilde{f}\eea
where
\bea\label{energy}
E = \lambda[\lambda+(2n-2)]~.\eea
Thus any holomorphic function $f$ of definite charge under $\xi$, or 
equivalently degree under $r\partial/\partial r$, 
corresponds to an eigenfunction of the Laplacian $-\nabla^2_L$ on the link. 
The charge $\lambda$ is then related simply to the eigenvalue $E$ by the above 
formula (\ref{energy}). We may thus in particular define the holomorphic spectral invariant
\bea
Z(t) = \sum_{i=0}^{\infty} \exp(-t\lambda_i)\eea
where $\{\lambda_i\}_{i=0}^{\infty}$ is the holomorphic spectrum, in the above sense. 
This is also the trace of a kernel, namely the \emph{Szeg\"o kernel}.

Recall the following classical result:
\begin{theorem}\label{weyl} (\cite{names}) Let $(L,g_L)$ be a compact Riemannian manifold with 
heat kernel trace $\Theta(t)$ given by (\ref{heat}). Then
\bea
\vol[g_L] = \lim_{t\searrow 0} \ (4\pi t)^{n-\tfrac{1}{2}} \ \Theta(t)~.
\eea
\end{theorem}
In the holomorphic setting we have:
\begin{theorem}\label{character}(\cite{MSY2}) For $g_L\in \mathcal{S}(X_0,\Omega)$ with Reeb vector field $\xi$ we have
\bea
\frac{\vol[g_L]}{\vol[S^{2n-1},g_{\mathrm{can}}]} = \lim_{t\searrow 0} \ t^n \ Z(t)~.\eea
\end{theorem}
This result first appeared for regular 
Sasaki-Einstein manifolds in \cite{herzog}. 
The proof of Theorem \ref{character} is essentially 
the Riemann-Roch Theorem. 
Suppose that $g_L\in\mathcal{S}(X_0,\Omega)$ is quasi-regular\footnote{The reader 
who is uneasy with orbifolds may take a regular Sasakian manifold in what follows. However,
the point here is that quasi-regular Reeb vector fields are dense in the space of 
Reeb vector fields, since the rationals are dense in the reals; regular Sasakian 
structures are considerably more special.}, so that $\xi$ 
generates a $U(1)$ action on $L$. The quotient is a \emph{Fano} orbifold $(V,\omega_V)$:
\begin{definition}
A compact K\"ahler orbifold $(V,\omega_V)$ is \emph{Fano} if the cohomology class 
of the Ricci-form in $H^{1,1}(V)$ is represented by a positive $(1,1)$-form. 
\end{definition} 
Holomorphic functions on $X_0$ that are eigenstates under $\mathcal{L}_{\xi}$ 
correspond to holomorphic sections of a holomorphic orbifold line 
bundle $\mathcal{L}^{-k}\rightarrow V$ for some $k$; here $\mathcal{L}$ 
is the associated holomorphic line orbibundle to the $U(1)$ principal orbibundle 
$U(1)\hookrightarrow L\rightarrow V$. This is holomorphic since the 
curvature is proportional to the K\"ahler form on $V$, which is of Hodge type $(1,1)$. The number of holomorphic sections 
is given by an orbifold version of the Riemann-Roch theorem, involving characteristic 
classes on $V$. In general this is rather 
more complicated than the smooth Riemman-Roch theorem, but the 
limit in Theorem \ref{character} simplifies the formula considerably: only 
a leading term contributes. 
The essential point now is that the volume may also be written in terms of 
Chern classes:
\begin{proposition}\label{quasi} For a quasi-regular Sasakian metric $g_L\in\mathcal{S}(X_0,\Omega)$ one has
\bea
\frac{\vol[g_L]}{\vol[S^{2n-1},g_{\mathrm{can}}]} = \frac{\beta}{n^n}\int_V c_1(V)^{n-1}~.\eea
Here $(S^{2n-1},g_{\mathrm{can}})$ is the round sphere metric; $\beta\in\mathbb{Q}$ is defined by
\bea
c_1(\mathcal{L})=-\frac{c_1(V)}{\beta}\in H^2_{\mathrm{orb}}(V;\mathbb{Z})\eea
where $H^2_{\mathrm{orb}}(V;\mathbb{Z})$ is the orbifold cohomology of Haefliger\footnote{One 
defines $H^*_{\mathrm{orb}}(V;\mathbb{Z}) = H^*(BV;\mathbb{Z})$ where $BV$ is the 
classifying space for $V$. For details, see for example \cite{boyer}.} \cite{haefliger}, 
$c_1(V)$ is the first Chern class of the holomorphic tangent bundle of $V$, and 
$\mathcal{L}$ is the orbifold line bundle associated to the $U(1)$ principal orbibundle 
$U(1)\hookrightarrow L\rightarrow V$. \end{proposition}
The cohomology group $H^2_{\mathrm{orb}}(V;\mathbb{Z})$ classifies 
orbifold line bundles over $V$, in exactly the same way that $H^2(V;\mathbb{Z})$ 
classifies line bundles when $V$ is smooth. The Proposition relates the volume of $g_L\in \mathcal{S}(X_0,\Omega)$, 
for quasi-regular $\xi$, to characteristic classes of $V$. Theorem \ref{character} then follows from the fact 
that quasi-regular Reeb vector fields are dense in the space of 
all Reeb vector fields, and vol is continuous.

We end this section with a localisation formula for the volume. We first 
require a 
\begin{definition}\label{resolve} Let $(X,\Omega)$ be a Gorenstein singularity with isolated singular 
point $p$, $X_0=\mathbb{R}_+\times L$, and let $\mathcal{S}(X_0,\Omega)$ be as above, with 
respect to an effective holomorphic action of $\mathbb{T}^s$. We say that an orbifold $\hat{X}$ is 
a \emph{partial resolution} of $X$ if 
\bea
\pi: \hat{X}\rightarrow X\eea
is a $\mathbb{T}^s$-equivariant map with $\pi:\hat{X}\setminus E 
\rightarrow X_0$ a $\mathbb{T}^s$-equivariant biholomorphism for some exceptional set $E$. If $\hat{X}$ is 
smooth we say that it is a resolution of $X$.\end{definition}
First note that such a partial resolution $\hat{X}$ always exists: one can take any quasi-regular 
Reeb vector field $\xi\in\mathtt{t}_s$ and blow up the corresponding orbifold $V$. 
In this case the exceptional set $E=V$ and the partial resolution is clearly equivariant. We then 
have the following
\begin{theorem}\label{localise}(\cite{MSY2}) Let $g_L\in \mathcal{S}(X_0,\Omega)$ and pick a partial resolution $\hat{X}$ of $X$.
Suppose that the K\"ahler form of $X_0$ extends to a suitable smooth family of 
K\"ahler forms on $\hat{X}$ (see \cite{MSY2} for details). Then
\bea
\frac{\vol[g_L]}{\vol[S^{2n-1},g_{\mathrm{can}}]} = \sum_{\{F\}} \frac{1}{d_F}\int_F 
\prod_{m=1}^R \frac{1}{\langle\xi,u_m\rangle^{n_m}}\left[\sum_{a\geq 0}\frac{c_a(\mathcal{E}_m)}{\langle\xi,u_m\rangle^a}\right]^{-1}~.\eea
Here
\begin{itemize}
\item $E\supset\{F\}=$ set of connected components of the fixed point set, 
where $\xi$ is a \emph{generic} vector $\xi\in\mathtt{t}_s$; that is, the 
orbits of $\xi$ are dense in the torus $\mathbb{T}^s$.
\item For fixed connected component $F$, the linearised $\mathbb{T}^s$ action 
on the normal bundle $\mathcal{E}$ of $F$ in $\hat{X}$ is determined by a set of 
weights $u_1,\ldots,u_R\in\mathbb{Q}^s\subset\mathtt{t}_s^*$. 
$\mathcal{E}$ then splits
$\mathcal{E}=\bigoplus_{m=1}^R \mathcal{E}_m$ 
where $\mathrm{rank}_{\mathbb{C}} \ \mathcal{E}_m=n_m$ and $\sum_{m=1}^R n_m = 
\mathrm{rank}_{\mathbb{C}}(\mathcal{E})$.
\item $c_a(\mathcal{E}_m)$ are the Chern classes of $\mathcal{E}_m$. 
\item When $\hat{X}$ has orbifold singularities, the normal fibre 
to a \emph{generic} point on $F$ is not a complex vector space, but 
rather an orbifold $\mathbb{C}^l/\Gamma$. Then $\mathcal{E}$ is more 
generally an orbibundle and $d_F=|\Gamma|$ denotes the order of $\Gamma$.
\end{itemize}
\end{theorem}
This theorem is proven as follows. One notes that the volume $\vol[g_L]$ may be 
written \cite{MSY2} as
\bea
\vol[g_L] = \frac{1}{2^{n-1}(n-1)!}\int_{X_0} e^{-r^2/2}\frac{\omega^n}{n!}~.\eea
The function $r^2/2$ is precisely the \emph{Hamiltonian function} for the Reeb vector field 
$\xi$. One may then naively apply the theorem of Duistermaat-Heckman \cite{DH1,DH2}, which easily 
extends to non-compact manifolds and orbifolds. This theorem 
localises such an integral to the fixed point set of $\xi$. However, 
since $\|\xi\|_g^2=r^2$, the integral formerly localises at the \emph{singular point} $p=\{r=0\}$ of $X$. 
To obtain a sensible answer, one must first resolve the singularity, as in Theorem \ref{localise}. 
Such a proof requires that the K\"ahler form of $X_0$ extends to a suitable smooth family of 
K\"ahler forms on $\hat{X}$. One then applies the Duistermaat-Heckman theorem 
to this family, and takes the cone limit $(X,g)$. The limit is independent of 
the choice of resolving family of K\"ahler metrics on $\hat{X}$.

\begin{remark} The author believes that the technical condition requiring the ability to 
extend the K\"ahler form $\omega$ on $X_0$ to a smooth family of K\"ahler forms on $\hat{X}$ is probably
redundant. In fact, one can also formerly apply the \emph{equivariant} Riemann-Roch 
theorem to $\hat{X}$, with respect to the holomorphic action of $\mathbb{T}^s$, and Proposition \ref{quasi} to obtain the same result. 
This has the advantage of not requiring existence of any K\"ahler metrics. \end{remark}

Note the theorem guarantees that $\vol:\mathcal{R}(X_0)\rightarrow\mathbb{R}_+$, relative to the volume of the round sphere,
is a rational function of $\xi$ with rational coefficients. Since vol is strictly convex, its critical 
points are isolated. By Theorem \ref{localise} one thus sees that the volume of a Sasaki-Einstein 
manifold, relative to that of the round sphere, is an \emph{algebraic number}.

Since the expression in Theorem \ref{localise} is
rather formidable, we end with an
\begin{example} Let $X$ be an affine toric variety with moment polytope  
$\mathcal{C}^*$. Let $\pi:\hat{X}\rightarrow X$ be a toric resolution of $X$, with 
some choice of K\"ahler metric, and denote the moment polytope by $\hat{\mu}(\hat{X}) = P\subset 
\mathtt{t}_n^*$. It is standard that such a resolution and K\"ahler metric always exist. Let $\mathrm{Vert}(P)$ 
denote the set of vertices of $P$. For each vertex $A\in \mathrm{Vert}(P)$ there are precisely 
$n$ outward-pointing edge vectors; these may be taken to be primitive vectors $u_A^i\in \mathbb{Z}^n\subset\mathtt{t}_n^*$, 
$i=1,\ldots,n$. This follows since the resolution is smooth. Then
\bea\label{toricvol}
\frac{\vol[g_L]}{\vol[S^{2n-1},g_{\mathrm{can}}]} = \sum_{A\in \mathrm{Vert}(P)} \prod_{i=1}^n \frac{1}{\langle \xi,u_A^i\rangle}~.
\eea
Note that this gives an unusual way of computing the Euclidean volume of the polytope $\Delta(\xi)$, which, by (\ref{euclid}), is 
also essentially the left-hand side of (\ref{toricvol}).
\end{example}

\section{Obstructions to the existence of Sasaki-Einstein metrics}\label{obstructions}

In this final section we examine the following 
\begin{problem}\label{prob} Let  $(X,\Omega)$ be a Gorenstein singularity with isolated singular 
point $p$, $X\setminus \{p\} = X_0 = \mathbb{R}_+\times L$, and let $g_L\in \mathcal{S}(X_0,\Omega)$ 
be a Sasakian metric with Reeb vector field satisfying (\ref{charge}). When does $X_0$ admit a Ricci-flat K\"ahler cone metric with 
this Reeb vector field?\end{problem}
We shall describe three obstructions. The first is a natural corollary of the previous section, whereas 
the remaining two 
obstructions are based on classical theorems in differential geometry, and in particular lead to new
obstructions to the existence of K\"ahler-Einstein orbifold metrics.

The previous section implies that a Sasaki-Einstein metric is a critical 
point of the volume functional, thought of as a function on the space $\mathcal{R}(X_0,\Omega)$ of Reeb vector fields 
satisfying (\ref{charge}). Thus the Reeb vector field in Problem \ref{prob} 
must be a critical point of this function, or equivalently of the Einstein-Hilbert action on $L$. 
We have already given this condition, in Proposition \ref{der}:
\bea\label{devagain}
\diff \vol (Y) = -n\int_L \eta(Y)~.\eea
Here $Y\in\mathtt{t}_s$ and recall that $\eta$ is the contact one-form on $L$. For a Sasaki-Einstein metric (\ref{devagain}) must be 
zero for all holomorphic vector fields $Y\in\mathtt{t}_s$ satisfying $\mathcal{L}_Y\Omega = 0$. 
Suppose that the Sasakian metric $g_L\in\mathcal{S}(X_0,\Omega)$ is quasi-regular. Then there is a quotient 
$V=L/U(1)$ which is a Fano orbifold. In this case we have the following
\begin{theorem}\label{futaki} Let $g_L\in \mathcal{S}(X_0,\Omega)$ be quasi-regular and let $Y\in\mathtt{t}_s$ 
with $\mathcal{L}_Y\Omega=0$. Then
\bea
\diff \vol(Y) = -\frac{\ell}{2} F \left(\mathcal{J}_V(Y_V)\right)~.\eea
Here $\ell=2\pi\beta/n$ is the length of the generic Reeb $S^1$ fibre, $\mathcal{J}_V$ is the complex 
struture tensor on $V$, $Y_V$ is the push-forward of $Y$ to $V$, and $F:\mathrm{aut}(V)\rightarrow \mathbb{R}$ 
is the \emph{Futaki invariant} of $V$.
\end{theorem}
The Futaki invariant is a well-known obstruction to the existence of a K\"ahler-Einstein metric on $V$ \cite{futaki}. 
It is conventionally defined as follows. Let $(V,\omega_V)$ be a compact K\"ahler orbifold with K\"ahler form 
satisfying $[\rho_V]=2n[\omega_V]\in H^{1,1}(V)$, where $\rho_V=\mathrm{Ric}(\mathcal{J}_V\cdot,\cdot)$
is the Ricci-form of $(V,\omega_V)$. Then there exists a real function $f$ on $V$, unique up to an additive constant, satisfying
\bea
\rho_V - 2n\omega_V = i\partial\bar{\partial}f~.\eea
Given any real holomorphic vector field $\zeta\in\mathrm{aut}(V)$ we define
\bea
F(\zeta) = \int_V \mathcal{L}_{\zeta}f \ \frac{\omega_V^{n-1}}{(n-1)!}~.\eea
Clearly, if $(V,\omega_V)$ is K\"ahler-Einstein of scalar curvature $4n(n-1)$ then $f$ is constant, 
and the function $F:\mathrm{aut}(V)\rightarrow\mathbb{R}$ vanishes. However, more generally $F$ also 
satisfies the following rather remarkable properties
\cite{futaki, calabi}:
\begin{itemize}
\item $F$ is independent of the choice of K\"ahler metric representing $[\omega_V]$.
\item The complexification $F_{\mathbb{C}}$ is a Lie algebra homomorphism $F_{\mathbb{C}}:\mathrm{aut}_{\mathbb{C}}(V)\rightarrow\mathbb{C}$.
\end{itemize}
Note that the first point is implied by Proposition \ref{reeby} and Theorem \ref{futaki}. 
The second is due to Calabi \cite{calabi}. The Sasakian setting thus 
gives a \emph{dynamical} interpretation of the Futaki invariant. 

We now turn to two further simple obstructions to the existence of solutions to Problem \ref{prob}. 
These are based on the classical theorems of Lichnerowicz \cite{lich} and Bishop \cite{bishop}, respectively. 
We begin with the latter. Recall from Proposition \ref{reeby} that the volume of a Sasakian 
metric $g_L\in \mathcal{S}(X_0,\Omega)$ is determined by its Reeb vector field $\xi$. 
In particular, one can compute this volume $\vol(\xi)$ using Theorem \ref{character}. In some cases 
one can compute the trace over holomorphic functions on $X$ directly. Now, Bishop's theorem \cite{bishop} implies 
that for any $(2n-1)$-dimensional Einstein manifold $(L,g_L)$ 
with $\mathrm{Ric}(g_L)=2(n-1)g_L$ we have
\bea
\vol(L,g_L)\leq \vol(S^{2n-1},g_{\mathrm{can}})~.\eea
Combining these two results we have
\begin{theorem}\label{bish}[Bishop obstruction] (\cite{GMSY}) 
Let $(X, \Omega)$ be as in Problem \ref{prob}, with Reeb vector field $\xi$. 
If $\vol(\xi)>\vol(S^{2n-1},g_{\mathrm{can}})$ then $X_0$ admits no 
Ricci--flat K\"ahler cone metric with Reeb vector field $\xi$.
\end{theorem}
{\it A priori}, it is not clear this condition can ever obstruct 
existence. We shall provide examples below. However, in \cite{GMSY} 
we conjectured that for \emph{regular} Reeb vector fields $\xi$ Theorem 
\ref{bish} never obstructs. This is equivalent to the following conjecture
about smooth Fano manifolds:
\begin{conjecture}\label{fano} (\cite{GMSY}) Let $V$ be a smooth Fano manifold 
of complex dimension $n-1$ with Fano index $I(V)\in\mathbb{N}$. Then
\bea
I(V)\int_V c_1(V)^{n-1} \leq n \int_{\mathbb{CP}^{n-1}} c_1(\mathbb{CP}^{n-1})^{n-1} = n^n
\eea
with equality if and only if $V=\mathbb{CP}^{n-1}$.
\end{conjecture}
This is related to, although slightly different from, a standard 
conjecture about Fano manifolds. For further details, see \cite{GMSY}.

We turn now to the Lichnerowicz obstruction. Suppose that $(L,g_L)$ is Einstein with $\mathrm{Ric}(g_L)=2(n-1)g_L$. The first non--zero 
eigenvalue $E_1>0$ of $-\nabla^2_L$ is bounded from below:
\bea
E_1 \geq 2n-1~.\eea
This is Lichnerowicz's theorem \cite{lich}. 
Moreover, equality holds if and only 
if $(L,g_L)$ is isometric to the round sphere $(S^{2n-1},g_{\mathrm{can}})$  \cite{obata}. 
From (\ref{energy}), we immediately see that for holomorphic functions $f$ on $X_0$ of 
charge $\lambda$ under $\xi$, Lichnerowicz's 
bound becomes $\lambda\geq 1$. This leads to a potential holomorphic obstruction to 
the existence of Sasaki--Einstein metrics:
\begin{theorem}\label{lich}[Lichnerowicz obstruction] (\cite{GMSY})
Let $(X, \Omega)$ be as in Problem \ref{prob}, with Reeb vector field $\xi$.
Suppose that $f$ is a holomorphic function on $X$ of positive charge 
$\lambda<1$ under $\xi$. Then $X_0$ admits no Ricci--flat K\"ahler cone metric with 
Reeb vector field $\xi$. 
\end{theorem}
Again, it is not immediately clear that this can ever obstruct existence.
Indeed, for \emph{regular} $\xi$ one can prove \cite{GMSY} that this never obstructs. 
This follows from the fact that $I(V)\leq n$ for any smooth Fano $V$ of complex 
dimension $n-1$. However, there exist plenty of obstructed quasi-regular examples.

Notice then that the Lichnerowicz obstruction involves 
holomorphic functions on $X$ of small charge with respect to $\xi$,
whereas the Bishop obstruction is a statement about the volume, which 
is determined by the asymptotic growth of holomorphic functions on $X$, analogously
to Weyl's asymptotic formula.

\begin{example} Our main set of examples of Theorems \ref{bish} and \ref{lich}
is provided by isolated quasi-homogeneous hypersurface singularities. 
Let $\mathbf{w}\in \mathbb{N}^{n+1}$ be a vector of positive weights. 
This defines an action of $\C^*$ on $\C^{n+1}$ via
\bea
(z_1,\ldots,z_{n+1})\mapsto (q^{w_1}z_1,\ldots,q^{w_{n+1}}z_{n+1})\eea
where $q\in\C^*$. Without loss of generality one can take 
the set $\{w_i\}$ of components of $\mathbf{w}$ to have no common factor. This ensures 
that the above $\C^*$ action is effective. Let
\bea
F:\C^{n+1}\rightarrow \C\eea
be a quasi--homogeneous polynomial on $\C^{n+1}$ with 
respect to $\mathbf{w}$. This means 
that $F$ has definite degree $d$ under the above $\C^*$ action:
\bea
F(q^{w_1}z_1,\ldots,q^{w_{n+1}}z_{n+1})=q^d 
F(z_1,\ldots,z_{n+1})~.\eea 
Moreover, we assume that the affine algebraic variety
\bea
X=\{F=0\}\subset \C^{n+1}\eea
is smooth everywhere except at the origin $(0,0,\ldots,0)$. 
For obvious reasons, such $X$ are called isolated 
quasi--homogeneous hypersurface singularities. 
The corresponding link $L$ is the
intersection of $X$ with the unit sphere in $\C^{n+1}$:
\bea
\sum_{i=1}^{n+1}|z_i|^2=1~.
\eea
We define a nowhere zero holomorphic $(n,0)$-form $\Omega$ 
on the smooth part of $X$ by
\be
\Omega = \frac{\diff z_1\wedge\cdots\wedge\diff z_n}{\partial F/\partial z_{n+1}}~.\ee
This defines $\Omega$ on the patch where $\partial F/\partial z_{n+1}\neq0$. One has 
similar expressions on patches where $\partial F/\partial z_i\neq0$ for each $i$, and it is simple to 
check that these glue together into a nowhere zero form $\Omega$ on $X_0$. 
Thus all such $X$ are Gorenstein, and moreover they come equipped 
with a holomorphic $\C^*$ action by construction. The orbit space of this 
$\C^*$ action, or 
equivalently the orbit space of $U(1)\subset \C^*$ on the link, is 
a complex orbifold $V$. In fact, $V$ is the weighted 
variety defined by $\{F=0\}$ in the weighted
projective space $\mathbb{WCP}^n_{[w_1,w_2,\ldots,w_{n+1}]}$. 
It is not difficult to show that $V$ is a Fano orbifold if and only if 
\be
|\mathbf{w}|-d>0\ee
where $|\mathbf{w}|=\sum_{i=1}^{n+1}w_i$. 
To see this, first notice that 
$|\mathbf{w}|-d$ is the charge of $\Omega$ under $U(1)\subset \C^*$. 
To be precise, if $\zeta$ denotes the holomorphic vector field on $X$ with
\be
\mathcal{L}_{\zeta}z_j = w_j i z_j\ee
for each $j=1,\ldots,n+1$, then
\be
\mathcal{L}_{\zeta}\Omega = (|\mathbf{w}|-d) i\Omega~.\ee
Positivity of this charge then implies \cite{MSY2} that 
the cohomology class of the natural Ricci--form induced on $V$ is represented by 
a positive $(1,1)$--form, which is the definition that $V$ is Fano. 
If there exists a Ricci--flat K\"ahler metric on $X_0$ which is a cone 
under $\R_+\subset \C^*$, then the correctly normalised 
Reeb vector field $\xi\in\mathcal{R}(X_0,\Omega)$ is thus
\bea
\xi = \frac{n}{|\mathbf{w}|-d}\zeta~.\eea

Bishop's theorem then requires, for existence of a Sasaki--Einstein 
metric on $L$ with Reeb vector field $\xi$, 
\bea\label{bishbound}
d\left(|\mathbf{w}|-d\right)^n\leq wn^n\eea
where $w=\prod_{i=1}^{n+1}w_i$ is the product of the weights. 
The computation of the volume that gives this inequality may be found in \cite{GMSY}. It is 
simple to write down infinitely many examples of isolated quasi-homogeneous hypersurface singularities 
that violate this inequality, and are thus 
obstructed by Theorem \ref{bish}. 

Lichnerowicz's theorem requires, on the other hand, that
\bea\label{lichbound}
|\mathbf{w}|-d\leq nw_{\mathrm{min}}
\eea
where $w_{\mathrm{min}}$ is the smallest weight. Moreover, this bound can be saturated if and only if $(X_0,g)$ is $\C^{n}\setminus\{0\}$ with 
its flat metric. It is again clearly trivial to construct many 
examples of isolated hypersurface singularities that violate this bound, and are hence 
obstructed by Theorem \ref{lich}. 
\end{example}

We conclude by making some remarks on the possible solution to Problem \ref{prob}. Let us first 
comment on the case of toric varieties. The Reeb vector field for a critical point of the Einstein-Hilbert action, 
considered as a function on $\mathcal{R}(X_0,\Omega)$,
exists and is unique, by Theorem \ref{torictheorem}. The remaining condition 
for a Ricci-flat K\"ahler cone metric may be 
written as a real Monge-Amp\`ere equation on the polytope $\mathcal{C}^*$ \cite{MSY}. This has recently been shown to always 
admit a solution in \cite{futakinew}. We state this as
\begin{theorem} (\cite{futakinew}) Let $(X,\Omega)$ be an affine toric Gorenstein singularity. Then $X_0$ 
admits a $\mathbb{T}^n$-invariant Ricci-flat K\"ahler cone metric, with Reeb 
vector field determined by Theorem \ref{torictheorem}.
\end{theorem}
Thus toric varieties are unobstructed. However, this still leaves us with the general non-toric case. 
Problem \ref{prob} is in fact closely related to a major open conjecture in the K\"ahler category. 
If $\xi$ is \emph{regular}, then existence of a Ricci-flat K\"ahler cone metric on $X_0$ with this Reeb vector field is 
equivalent to existence of a K\"ahler-Einstein metric on the Fano $V$, and we then have 
the following conjecture due to Yau:
\begin{conjecture} (\cite{yau}) A Fano manifold $V$ admits a K\"ahler-Einstein metric 
if and only if it is \emph{stable} in the sense of Geometric Invariant Theory.
\end{conjecture}
A considerable amount of progress has been made on this conjecture, notably by 
Tian and Donaldson. However, it is still open in general. The conjecture 
is closely related to Problem \ref{prob}, and it is clearly of interest to 
extend the work in the K\"ahler setting to the Sasakian setting. 

Physics
might also provide a different viewpoint. Suppose that $X$ admits a \emph{crepant resolution} 
$\pi:\hat{X}\rightarrow X$. This means that the holomorphic volume form 
$\Omega$ on $X_0$ extends smoothly\footnote{This is certainly 
an extra constraint. For example, the singularities $w^{2k+1}+x^2+y^2+z^2=0$ admit  
no crepant resolution by an argument in \cite{GMSY}.} as a holomorphic volume form onto the resolution $\hat{X}$. 
Then one expects the derived category of coherent sheaves $\mathcal{D}^b(\mathrm{coh}(\hat{X}))$ on $\hat{X}$ 
to be equivalent to the derived category of representations $\mathcal{D}^b(\mathrm{Reps} \ (\mathcal{Q},\mathcal{R}))$ of a quiver 
$\mathcal{Q}$ with relations $\mathcal{R}$. In fact this could be stated formerly as a conjecture; it is known to be true for various sets of examples. A quiver $\mathcal{Q}$ is simply a directed graph. $\mathcal{R}$ is a set of relations on the path algebra $\C\mathcal{Q}$ of the quiver. 
The above correspondence between derived categories, if correct, would allow for a more precise 
mathematical statement of what the AdS/CFT map is. Physically, one is placing 
D3-branes at the singular point $p$ of $X$; mathematically, a D3-brane at a point on 
$X$ corresponds to the structure sheaf of that point. The corresponding 
representation of $\mathcal{Q}$ then defines 
an $\mathcal{N}=1$ supersymmetric quantum field theory in four dimensions; this is 
precisely the quantum field theory on the D3-brane. 
The quiver representation determines the gauge group and matter content of the quantum field theory, while the relations specify 
the superpotential, which determines the interactions. The choice of $(\mathcal{Q},\mathcal{R})$ is known 
to be non-unique, and this non-uniqueness is related to a duality known as Seiberg duality \cite{seiberg}. 
The AdS/CFT correspondence then implies that $X$ admits a Ricci-flat 
K\"ahler cone metric (in dimension $n=3$) if and only if this supersymmetric quantum 
field theory flows to a dual infra-red fixed point under renormalisation group flow; this infra-red 
fixed point is precisely the superconformal field theory in the AdS/CFT correspondence.
Thus the solution to Problem \ref{prob}, in complex dimension $n=3$, is related to the low-energy behaviour of certain 
\emph{supersymmetric quiver gauge theories} in four dimensions.


\subsection*{Acknowledgments}
\noindent I would very much like to thank the organisers of the conference for inviting me to 
speak. I also thank my collaborators: J.~Gauntlett, D.~Martelli, D.~Waldram and S.-T.~Yau.
I am supported by NSF grants DMS--0244464, DMS--0074329 and DMS--9803347.

\end{document}